\documentclass[12pt]{article}
\usepackage{amsmath,amsthm,amsfonts}

\def\be{\begin{equation}}
\def\ee{\end{equation}}

\def\R{{\sf I\kern-.2em R}}
\def\N{{\sf I\kern-.2em N}}
\def\C{\kern.1em{\raise.47ex\hbox{$\scriptscriptstyle
$}}\kern-.40em{\sf C}}
\def\Z{{\sf Z\kern-.32em Z}}

\def\be{\begin{equation}}
\def\ee{\end{equation}}

\def\limind{\mathop{\rm lim\,ind}_{\rightarrow}}
\def\limproj{\mathop{\rm lim\,proj}_{\leftarrow}}

\newtheorem{theorem}{\noindent Theorem}
\newtheorem{lemma}{\noindent Lemma}
\newtheorem{definition}{\noindent Definition}

\author{A.~M.~Vershik}
\title{Gelfand--Tsetlin algebras, expectations,
inverse limits, Fourier analysis}
\date{}

\begin{document}

\maketitle

\abstract{This text mainly follows my talk at the
 conference ``Unity of Mathematics'' devoted to the 90th anniversary
 of I.~M.~Gelfand (Harvard, September 2003). I introduce
 some new notions that are related to several old ideas of
 I.~M.\ and try to give a draft of the future development of this area,
 which includes the representation theory of inductive families of
 groups and algebras and Fourier analysis on such groups.
 I also include a few reminiscences about I.~M.\ as my guide.}

\section*{0\ \ \  Historical excursus: I.~M.~Gelfand as my
correspondence advisor}

The first substantial series of mathematical works that I had
studied being a student was the series of papers
by Gelfand, Raikov, and Shilov (GIMDARGESH, as I called it to myself)
on commutative normed rings and subsequent papers on the generalized
Fourier analysis. This theory became a mathematical shock for me,
I was striked by its beauty and naturalness, universality
and depth.

Before this I hesitated whether I should join
the Chair of Algebra --- I attended the course of Z.~I.~Borevich
on group theory and the course of D.~K.~Faddeev on Galois theory ---
or the Chair of Mathematical Analysis, where my first advisor
G.~P.~Akilov worked; in the latter case I could choose the complex
analysis (V.~I.~Smirnov, N.~A.~Lebedev) or the real and
functional analysis (G.~M.~Fikhtengolts, L.~V.~Kantorovich).
But now the choice was clear: the functional analysis.
At the same time I was more interested in the works of the Moscow
(Gelfand's) school of functional analysis focused on noncommutative
problems than in the works of the Leningrad school that was oriented rather
towards the theory of functions and operator theory.

Since then the works of I.~M.~Gelfand and his school in various
fields became a kind of mathematical guidebook for me. My master
thesis was devoted to the theory of generalized functions; this
topic equally interested the Leningrad mathematicians
(L.~V.~Kantorovich, G.~P.~Akilov). Later, following G.~P.~Akilov's
advice, I began to study the representation theory, which was
at the time absolutely unrepresented in Leningrad. While
my being a postgraduate student, I.~M.~Gelfand popularized
problems concerning measure theory in infinite-dimensional
spaces, inspired by the theory of distributions,  the notion
of generalized random processes, and quantum physics.
These problems were communicated to us by D.~A.~Raikov,
who, following I.~M.~Gelfand's advice, worked at the
new theory of locally convex and nuclear spaces,
which we also studied in G.~P.~Akilov's seminar.
In Leningrad, the measure theory
in linear topological spaces was being studied in
the late 50s -- early 60s
 by V.~N.~Sudakov and me.
At the time everybody believed that the theory of generalized functions and
measure theory in infinite-dimensional spaces would require
to overstep the limits of the conventional Banach functional analysis, which
would be replaced by the theory of nuclear spaces (Minlos--Sazonov
and Gelfand--Kostyuchenko theorems, quasi-invariant measures, etc.).
However, soon it became clear that the measure theory
in linear spaces is a natural part of the general
measure theory, and the Banach analysis continued to
be the traditional language of functional analysis.
After all, the interest to all these problems gradually
died away.

V.~A.~Rokhlin's arrival at Leningrad thoroughly changed the mathematical
landscape at the Department of Mathematics. In particular, he organized
seminars on ergodic theory and topology. V.~A.\ became my principal
advisor during my postgraduate studies and several subsequent years.
I seriously studied the theory of dynamical systems and general
measure theory, and both my dissertations were devoted to these problems.
But representation theory continued to fascinate me equally.
Even earlier, in his talks on problems of functional analysis
at the All-Union Conference on Functional Analysis and the 3rd
Mathematical Congress (1956), I.~M.~spoke about the von Neumann factors
and Wiener measure as about subjects that were
possibly related  and underestimated at the time. Later, in the 60s,
I began to study factors and relations of the theory of $C^*$-algebras,
introduced by Gelfand and Naimark, with the theory of dynamical
systems; this became the subject of my research for several years.

But for several short discussions with I.~M.\ in the mid and late 60s
and the correspondence acquaintance via V.~A.~Rokhlin (and possibly
via Yu.~V.~Linnik),
our close acquaintance took place
in the spring 1972. After a session of his seminar,
I began to talk to him about my work (joint with A.~A.~Shmidt)
on the limit statistics of cycles of random permutations;
and the next day, at his home, about my plans to study the
representations of the symmetric groups. Though first he said that
with them everything is clear
and started to talk enthusiastically
about the theory of symmetric functions, but later he agreed that
not everything is that clear and advised me to look at the paper
by E.~Thoma on the characters of the infinite
symmetric group, which was of most interest for me.
This paper played an important role in our subsequent studies
of this group with
my pupil S.~V.~Kerov. One of our principal contributions
was an explanation and a new proof of Thoma's result in terms of
representation theory (asymptotics of Young diagrams).
And in that conversation I.~M.\ approved
wholeheartedly of my ideas, which I later
called the asymptotic representation theory;
and even when he retold them to D.~Kazhdan, who
had appeared a little later,
he referred to the theorems on the asymptotic behavior of
Young diagrams, charaters, etc., which were only
conjectured at the time (many of them were proved
later in joint papers with S.~V.~Kerov), as to results already obtained.
Those of the results I talked about that were already proved
related  rather to probability theory (Poisson--Dirichlet
measures) and theory of dynamical systems than to representation theory.
Other groups besides the symmetric groups and their representations
were not discussed in those conversations.
I took leave of him and was going to depart for Leningrad.

Suddenly, at the day of my departure, I.~M., having found out, in a rather
complicated way, the phone number of my friends with whom I stayed at Moscow,
called me and asked to come to him immediately. He also invited
M.~I.~Graev, and during our long walk told me about the problem
of constructing the noncommutative integral of representations
for semisimple groups, and especially for $\rm{SL}(2,\mathbb R)$.
He had earlier offered this problem to other his pupils, but he said that he had no
doubt that it ``fitted'' me. I was slightly surprised, because
I supposed that he could not know to what extent I was acquainted with
the representation theory of Lie groups, and in particular
that of $\rm{SL}(2)$; as I have mentioned above, we did not discuss
these matters at all. But I.~M.\ was right --- this problem was offered to me at
a very appropriate moment. Several years before this conversation, at the youth seminar
organized by L.~D.~Faddeev and me, we studied Gelfand's volumes
on generalized functions and other useful things, which were not
widespread in Leningrad. And in the early 70s, apart from
my studies of the ergodic theory, I gave a course and seminars
just on the representation theory of groups and algebras,
tensor products, and factors.
Apparently, I.~M.\ had heard about it, but I did not ask him. Thus
his problem appeared at an appropriate moment. We coped with it within
several months (the end 1972 --- the beginning 1973).
The first paper in ``Uspekhi'' (``Russian Math. Surveys'') appeared
in a volume dedicated to Kolmogorov in 1973, and this was the
beginning of our collaboration with I.~M.\ and M.~I.~Graev, which
lasted with intervals about ten years and which I am going to
describe one day in more detail. That first (the best, in I.~M.'s
and my opinion) paper of this series was devoted to the ``integral''
of representations of $\rm{SL}(2,\mathbb R)$ and touched upon many
topics that are still actual; in that paper we rediscovered
several constructions that had recently appeared (Araki's Gaussian construction,
cohomologies in groups without Kazhdan's property, etc.), gave
the first explicit formulas for nonzero cohomologies of semisimple
group of rank 1, and constructed irreducible nonlocal
representations of current groups with values
in finite-dimensional Lie groups.
I.~M.\ repeatedly
(and the last time --- at this conference (Harvard 2003))
expressed his wish to continue our joint work in
this direction. We had no doubt that this series of papers
has various applications, which has already been repeatedly confirmed,
and the work would be continued.

\bigskip

This paper is devoted to a subject from another line, which
also goes back to I.~M.'s works. Having worked for many years with
inductive families of semisimple algebras, S.~V.~Kerov and the author
at once appreciated the importance of the notion that we called
the Gelfand--Tsetlin algebras; this notion is a generalization of the
well-known and still popular construction of the Gelfand--Tsetlin
bases for the unitary and orthogonal groups. These algebras
serve as a basis for the harmonic analysis and Fourier analysis on
noncommutative groups. They play an especially important role in
the representation theory of locally finite groups, symmetric groups,
and, more generally, inductive limits of groups and algebras.
Our joint works with A.~Okounkov (see
\cite{OV} and \cite{VO}) show how applying these algebras, and especially
a natural basis in them (the Young--Jucys--Murphy basis), allows one
to reconstruct the represenation theory of the symmetric
groups at a completely dufferent basis.
In my talk and in this paper I draw attention to yet another
idea, closely related to the previous one; namely, to the idea
of inverse limits of algebras with respect to conditional
expectations. For the symmetric group, this question
will be considered in detail in the joint work with N.~V.~Tsilevich
(under preparation). On the other hand, inverse limits of finite-dimensional
algebras generalize von Neumann's theory of complete and noncomplete
tensor products \cite{vN}, and I remember one of my first visits
to Gelfand's seminar in the late 50s, when
this von Neumann's paper was being discussed and commented by the head
of the seminar. In this paper I do not touch upon
one subject that I mentioned in the talk, namely, the results on
representations of the group of infinite matrices over a finite
field, which we intensively studied with S.~V.~Kerov during several last
years. It will be considered in other publications under preparation.

\section{Definition of a generalized expectation on a subalgebra}

  Let $A$ be a  $C^*$-algebra over {\bf C} with
  involution $*$,  and let $B$ be its involution $C^*$-subalgebra.
  All algebras in the paper are supposed to be algebras with identity, and
  all subalgebras are supposed to contain this identity.
  Here we mainly consider finite-dimensional algebras,
  but the definitions below are valid for the general case.

\begin{definition}
  A linear operator $P:A \longrightarrow B$ is called
  a {\em conditional mathematical expectation},
  or {\em expectation}\footnote{The word
  ``conditional'' is the traditional one, but  I prefer
  to omit it below, as well as the word
  ``mathematical,'' violating the old tradition. The reason
  is that the ``unconditional expectation'' is simply
  the ``conditional expectation'' onto the algebra of scalars
  $\bf C$ (conditions are trivial), thus if we fix a subalgebra
  $B$, we do not need to use the word ``conditional,'' because it
  is clear what are the ``conditions.''}
  for short, of the algebra $A$ onto the subalgebra $B$
  if
  \begin{enumerate}
\item  $P(b)=b$ and $P(b_1 a b_2)=b_1 P(a) b_2$
  for all $a \in A$ and $b,b_1,b_2 \in B$.

\item  $P(a^*)= a^*$, $P1=1$,

  and

\item   $P \ge 0$, which means that
  for all $a\in A$, $P(aa^*)$ is positive, i.e.,
  belongs to the real cone in $B$ generated by
  the elements of the form $bb^*$.
\end{enumerate}
We will say that $P$ is a {\em generalized expectation} if only
  the first and second conditions hold, and $P$ is a {\em true expectation}, or
{\em  expectation}, if condition 3 also holds.
\end{definition}

    The notion of (``conditional''!) expectation is well known
    and has been used in
  many situations; for commutative algebras, it coincides with the ordinary
  notion of (mathematical) conditional expectation on sigma-subfield or subalgebra.
  A fruitful example of {\it generalized}, i.e., {\it nonpositive} expectation
  appeared, I believe, only recently, in the very concrete situation of
  the group algebra
  of the symmetric group (see below), and this is the reason for considering
  this notion in full generality. Sometimes people require that an expectation
  $P$ should be not only positive  but even {\it totally positive}, but we will
  not put emphasis on this.

  Note also that it is clear from definition that the set of all expectations
  in an algebra $A$ to a subalgebra $B$ is a convex set.

  In the main part of the paper, our attention will be  focused
  not on a  single generalized expectation for some pair $B \subset A$, but
  on {\it sequences of generalized expectations in an
  inductive family of algebras}.

  It is not difficult to describe all expectations for finite-dimensional
  semisimple $C^*$-algebras over $\bf C$, which are the sums of several copies of
  full matrix algebras $M_n(\bf C)$, as well as to describe generalized
  conditional expectations for these algebras.
   Recall that for a general pair $(A,B)$,  where
   $A=\sum_{j=1}^m A_j$ is a finite-dimensional $C^*$-algebra,
   $B=\sum_{i=1}^k B_i$ is its $C^*$-subalgebra, and
   $A_j= M_{k_j}({\bf C})$, $j=1, \dots, m$,
   $B_i=M_{n_i}({\bf C})$, $i=1, \dots, k$, are their
   decompositions into simple
   algebras, one can define a {\it bipartite multigraph} in which
  the first (upper) part of vertices is indexed by the subalgebras
  $B_i$, $i=1, \dots, k$,
  and the second (lower) part of vertices
  is indexed by the subalgebras $A_j$, $j=1, \dots ,m$,
  and the multiplicity of an edge $(i,j)$ is equal to the
  number of copies of the subalgebra $B_i$ as a
  subalgebra of $A_j$.
  We will use this construction in the below theorem (claim 2).
  For the sake of clarity, we consider the multiplicity-free case when each $B_i$ belongs
  to at most one $A_j$; a pair $(i,j)$ is called admissible if it is an
  edge, or $B_i \subset A_j$. In order to determine the pair $(A,B)$
  uniquely  up to isomorphism, we must fix this bipartite
  multigraph and positive integers in each upper vertex
  (the dimensions of $B_i$).

  \begin{theorem}
  {\rm1.} First suppose that
  $A=M_n(\bf C)$ and its subalgebra $B$ is also
  a full matrix algebra $B=M_m(\bf C)$
  (that is, the multigraph reduces to two vertices and one edge).
  Then there exists a unique expectation
  $P(a)=pap$, where $a \in A$ and $p$ is the natural orthogonal projection
  determined by the identity of the algebra $B$.

   {\rm 2.} Suppose that $A$ is a finite-dimensional semisimple algebra
   and $B$ is its semisimple
    subalgebra as above. Then every conditional expectation
    $P:A \longrightarrow B$ is the
   sum $$P=\sum_{i,j} P_{i,j}$$ over all admissible pairs $(i,j)$
  of generalized expectations from claim {\rm1}:
  $P_{i,j}: A_j \longrightarrow B_i$, $P_{i,j}(a)=\lambda_{i,j} p_{i,j}ap_{i,j}$,
   where $\lambda_{i,j}$ are real numbers (for a true expectation,
   nonnegative real numbers)
   such that
   $\sum_j \lambda_{i,j} =1$ for every $i$.
\end{theorem}

   The proof of claim 1 is obvious; in order to prove claim 2,
   it suffices to separate  the restrictions of $P$ to each $A_j$
   by linearity of
   $P$ and then apply claim 1 and
   condition 2 from the definition of
   expectation ($P1=1$).

  Thus a real matrix $\{\lambda_{i,j}\}$ satisfying the condition
  $\sum_j \lambda_{i,j} =1$ for every $i$
  is a parameter on the set of generalized conditional expectations
  for a fixed semisimple finite-dimensional algebra and its subalgebra;
  for true expectations, we have an additional condition $\lambda_{i,j} \ge 0$,
  and $\{\lambda_{i,j}\}$  is a {\it Markovian matrix on the bipartite graph}.
  For this reason, in the general case
  we will say that the matrix $\lambda_{i,j}$
  is a generalized Markovian matrix. It is clear that the set of
  (generalized) expectations for a finite-dimensional
   pair $B\subset A$ is always nonempty.

 The conjugate operator to a generalized expectation $P$ is an operator
 $P^*$ from
 the space $A^*$ conjugate to $A$ to $B^*$. If $P$ is a true
 (positive) expectation, then $P^*$ maps each state (= positive normalized
 functional) on $B$ to some state on $A$. But since $P$ is
 not a homomorphism of algebras, it does not map traces
 (characters) to traces. We may consider more refined properties
 of expectations in regard to this fact, e.g., call an expectation
 {\it central} if the image of each trace is a trace, etc.
 We will not discuss this topic here.

The following natural question arises.
Suppose that
$P_B$ is an expectation for a pair of finite-dimensional
algebras
$A,B$. Let us regard $A$ as a vector space. The problem is to describe
the $*$-algebra
$E=<A,P_B>$ generated by the left action of $A$ and
$P_B$ in
${\rm END}(A)$. We give the answer
to this question in terms of the
decomposition of $E$
into simple algebras.

\begin{theorem}
Let  $\Gamma (L_B,L_A)$ be the bipartite graph
corresponding to the
pair $(A,B)$, where $L_A$ ($L_B$) are the vertices of $\Gamma$
corresponding
to the decomposition of $A$ ($B$), respectively. Then
the diagram of the
triple of algebras $(B \subset A \subset E)$ is the graph
$\Gamma
(L_B,L_A,L_E)$, where the bipartite part $\Gamma(L_A,L_E)$
is the {\it
reflection} of  $\gamma (L_B,L_A)$, which means that
$L_D \equiv
L_B$ and the edges between the vertices of $(L_A,L_E)$ are
the same as the
corresponding edges in $\Gamma (L_B,L_A)$. This means, in
particular, that the
algebra $E=<A,P_B>$ does not depend on the choice of the
expectation $P_B$,
but only on the subalgebra $B$ itself, so that we can denote it
by $E(A,B)$.
\end{theorem}

The proof of this theorem uses Theorem 1 (the structure of
expectations); we will
not give examples and details here.

\section{Two classes of examples for group algebras}

For the group algebras (over $\bf C$) of finite groups, we
present two types of expectations related to the group
structure. Since a linear map in the group algebra is determined by
its values on the group, we can state the question in terms of the
group.

\smallskip
1. The first type of examples relates to the case when
the value of the expectation at an element of the group (regarded as a subset of the group algebra)
belongs to the group again.

In this case we can formulate a pure group-theoretical question
concerning a group analog of expectation.

 Assume that $G$ is a finite group and $H$ is its subgroup.
 When there exists a map $p$ from $G$ onto $H$ such that
 $$p(h)=h, \qquad p(h_1gh_2)=h_1p(g)h_2,\quad p(e_G)=e_H$$
 for all $h, h_1, h_2 \in H$ and $g \in G$, where $e_H$ and $e_G$ are
 the identity elements in $G$ and $H$, respectively?

If such a map $p$ exists, we
say that it is a {\it virtual projection} of
 the group $G$ to the subgroup $H$.

 \begin{theorem}
The following two conditions are equivalent:
\begin{enumerate}
\item There exists a virtual projection $p:G \to H$.

\item There exists a set $K \subset G$ such that
\begin{itemize}
  \item[\rm a)] $K$ is invariant under the inner automorphisms
  generated by the elements of $H$, that is, for every
  $h\in H$, for every $k\in K$, $hkh^{-1}\in K$;

  \item[\rm b)] the intersection of the set $K$ with any left (equivalently, right)
   coset of $H$ in $G$ has only one element; in other words, for
  all $k,k' \in K$, $k\ne k'$, we have $k^{-1}k' \notin H$.
\end{itemize}
\end{enumerate}
\end{theorem}

\begin{proof}
   {The proof is straightforward, and we only supplement it with some comments.
Condition b) means that the group
$G$ can be partitioned into left cosets of the subgroup $H$, each
of them containing exactly one element of the set $K$;
thus $G\cong H\times K$, and for every $g \in G$
there is a unique left decomposition $g=hk$ with $h\in H$, $k\in K$;
condition a) gives
the right decomposition with the same $h$ but another
   $k' \in K: g=k'h$.
   We assert that there is a bijection between the set of all virtual
   projections $p:G\to H$ and the set of all subsets $K$ that satisfy
   these conditions.
   Namely, if $K$ enjoys properties a), b) above, then
   the corresponding virtual projection $p$ is given by the formula
         $$p(g)=h$$
     for the element $g=hk=k'h$;
   and vice versa: if $p$ is a virtual projection, then the
   set $K=p^{-1}(e_H) \subset G$ enjoys properties a), b).}
   \end{proof}

\smallskip\noindent{\bf Remark 1.}
It is clear from construction that the set $K$ is the union
of orbits of the group of inner automorphisms of $H$. If $O$ is
one of such orbits in $K$, then its characteristic function
commutes with $H$.
In the case of the symmetric group, $K$ is a single orbit.

\smallskip\noindent{\bf Remark 2.}
The set $K$ above can be described in the following terms
(E.~Vin\-berg's observation):  $\bar K=\{k \in G:\, k$ belongs to the
center of the group $H \cap k^{-1}Hk\}$.
Then our $K$ is a subset of $\bar K$, which is
$H$-invariant and intersects each left (and, automatically,
right) coset of the subgroup $H$ at one point.
\smallskip

For different groups, it may happen that such a set $K$ either is
nonunique, or does not exist at all.

It is an interesting question for what pairs $H \subset G$ a virtual
 projection (in terms of Theorem 3, a set $K$ with properties a), b))
 does exist. In the trivial example $G$ is the
 direct product of two groups: $G=H \times K$.

 As a nontrivial example, consider the symmetric groups
 $G=S_n$ and $H=S_{n-1}$ with the ordinary embedding; then $K$ is the set
 of transpositions $(i,n)$, where $i$ runs over $1,2, \dots, n$.
 The map $p:G \to H$ determined by this decomposition
 is a virtual projection; it simply deletes the
 element $n$ from a permutation. This projection was defined in \cite{KOV} (see
 also \cite{KOV1}) and called the {\it virtual projection}. It is easy to check
 that for $n>4$, the virtual projection and the corresponding set $K$ are
 unique; for $n=3,4$,  there are several possibilities to choose such a
 set $K$.

 Let us extend a virtual projection by linearity to an operator $P$ in the group algebra:
 $$P:C(G)\longrightarrow C(H)\subset C(G).$$

\begin{lemma}{The linear operator $P$ defined above is a generalized expectation
of the algebra $C(G)$ to $C(H)$ in the sense of Definition 1}.
\end{lemma}

An important remark: in general,
the generalized expectation $P$ does not satisfy the
positivity condition 3 from Definition 1;
for example, in the case of the symmetric group (see above),
this operator is not positive,
because, e.g.,
the signature of a permutation can change under this projection.
Thus $P$ is not an expectation, but a generalized expectation.

Thus we have defined a particular class of generalized
expectations on group algebras,
which arise from virtual projections on groups. A very
interesting problem is to describe pairs $(G,H)$ for which a
virtual projection, and hence the corresponding generalized expectation
on the group algebra, does exist. For an abelian group, it is
easy to describe all virtual projections (they exist for all pairs
$(G,H)$ and determine true expectations), but even for metabelian
groups, I do not know the answer.

For some classes of groups, such as free groups, ``local groups'' (see
\cite{V1}), Coxeter groups, presumably the following
recipe works: suppose that $G_n=<\sigma_1, \sigma_2, \dots, \sigma_n>$
and $G_n\supset G_{n-1}=<\sigma_1, \sigma_2, \dots, \sigma_{n-1}>$.
There exists a normal form of each element
of $G_n$ as a word in the alphabet $\sigma_1,\dots,
\sigma_n$ such that the deletion of the letter $\sigma_n$
in this normal form is a virtual projection
of $G_n$ onto $G_{n-1}$. This is true for free, locally
free, and symmetric groups (such a normal form does
exist).

\smallskip
2. The second type of examples is closer to the classical
definitions, because it leads to true (positive) expectations. Again
let $G$ and $H$ be a finite group and its subgroup,
respectively; now we allow the values of expectations at the
elements of the group not only to belong to
the group, but also to be equal to zero.
Define a projection
 $$P:C(G)\longrightarrow C(H)\subset C(G)$$
as follows: $P$ is the linear extension to the whole group algebra of the
following map on the group:  $P(h)=h$ for all $h \in H$, and
$P(g)=0$ if $g\in G$, $g \notin H$. This definition makes
sense for an arbitrary group and its subgroup.  Obviously, $P$ is a
(positive) expectation. For some reason, we
call it the  {\it
Plancherel expectation}. This definition leads, in particular, to
the Fourier analysis on the symmetric groups, which will be the subject
of the joint paper with N.~Tsilevich, which is now in preparation.

  It is easy to formulate the analog of Lemma 2
for algebras: the set of all generalized expectations
  $P:A \longrightarrow B$ is in a one-to-one correspondence with
the set of subspaces $T$ of $A$ satisfying the following properties:
\begin{enumerate}
\item $T$ is a closed complement to the subspace $B$ of the vector
space $A$;

\item $BTB \subset T$.
\end{enumerate}
 The correspondence is as follows: $T=\ker P$.

Because of the convexity of the set of expectations, we can consider
convex combinations of these two types of examples. For the
symmetric group, such deformations are related to the content of the papers
\cite{KOV, KOV1}.

\section{Gelfand--Tsetlin (GZ-) algebras}

 Now we introduce the central notion of the theory of
 inductive families of algebras (not only finite-dimensional).
 This notion follows the idea of the classical papers
 by Gelfand and Tsetlin \cite{GZ1,GZ2}, in which a particular basis
 was defined for the orthogonal $SO(n)$ and unitary $SU(n)$
 groups. This basis appears only if we consider
 not just one group, say $SO(n)$ or $SU(n)$, but
 the whole inductive family
 $SO(2) \subset SO(3)\subset \dots \subset SO(n)$ or
 $SU(1)\subset SU(2)\subset \dots  \subset SU(n)$ simultaneously.
 Since the restrictions
 of irreducible representations of the group $SO(n)$ to the subgroup
 $SO(n-1)$ (and similarly with $SU$)
 are multiplicity-free,
   this inductive family
 determines a basis  (Gelfand--Tsetlin basis), which
 is unique up to scalar multipliers
(see below). But even more important is the notion of
Gelfand--Tsetlin algebras, which was introduced for a general inductive
 family of algebras in our papers with S.~Kerov (a detailed exposition
is given in
\cite{VK}) and independently, but not in the same spirit, in
\cite{SV}. I do not know any papers about the Gelfand--Tsetlin algebras
even in the classical case (that of the universal enveloping algebras of
semisimple Lie algebras) apart from the paper \cite{Vin},
which concerns a  completely different problem.
 The most important problem is to define reasonable multiplicative
 generators of the Gelfand--Tsetlin algebras in terms of the initial
 algebras; having such generators, one can create
 the representation theory of the inductive family of algebras
 in a very natural way. The realization of this plan allows one to define an
 analog of the Fourier transform
 for algebras with inductive family of subalgebras inside it.
 For the symmetric group, these generators were
 defined (independently of GZ-algebras) by A.~Young and in more recent times
 by Jucys and Murphy (YJM-generators). The consistent development
 of the representation theory of the symmetric groups
 was given in \cite{OV, VO}. For other groups (even for the orthogonal
 and unitary groups), this is still not done. Below we consider only
 complex $*$-representations of algebras over $\bf C$.

\begin{definition}[Gelfand--Tsetlin algebra]
{Suppose we are given a finite or infinite family $A_k$, $k=0,\dots, n$
(here $n$ can be
finite or infinite), of semisimple algebras
over $\bf C$, $A_0=\bf C$, $A_k \subset A_{k+1}$.
Assume for the sake of clarity that the multiplicity of the restriction of
an irreducible representation of $A_k$ to $A_{k-1}$ for $k=1, \dots,
n-1$ is equal to one or zero (the so-called simple spectrum). By definition,
the Gelfand--Tsetlin
algebra $GZ_n$ is the algebra generated by the
{\it centers}, which we denote by $\zeta (A_k)\subset A_k$, $k=0,
\dots, n$:
 $$GZ_n=<\zeta (A_1), \dots, \zeta (A_n)>$$
 (the notation $<...>$ stands for the subalgebra of $A_n$
 generated by the contents of
 the brackets).}
\end{definition}

It is clear from this definition that all $GZ_k$ are abelian algebras
and the family of algebras $\{GZ_k\}_1^n$ is an inductive
family of subalgebras  in $A_n$ (the centers do not
form an inductive family);
the definition and the assumption on the simplicity of the spectrum
also imply that $GZ_n$ is a maximal abelian subalgebra of $A_n$.
Moreover, from the definition we can conclude that there is a particular
basis (defined up to scalars) in the algebra $GZ_n$, which we call
the GZ-basis; and, consequently, there is a particular basis in each irreducible
representation of $A_n$ --- this is what people usually called
the Gelfand--Tsetlin basis. In the case of the groups $SO(n)$ and $SU(n)$, this is
just the classical Gelfand--Tsetlin basis \cite{GZ1,GZ2}. It leads to
the well-known notion of Gelfand--Tsetlin patterns.

The elements of the GZ-basis of the algebra $GZ_n$ in the general case
are defined as such elements that each of them has a nonzero
image in only one irreducible representation. All
such elements are defined uniquely (up to scalar). We may say that
there is a bijection between this basis and paths in the graph of
the Bratteli diagram of the algebra $A_n$ (see below). As we have mentioned
above, a nontrivial problem is to describe the $GZ$-algebra, as well
as the $GZ$-basis, using some multiplicative generators of $GZ(A_n)$,
not in terms of representations, but in
intrinsic terms of the initial definition of the algebras $A_n$ (or groups
in the case when $A_n$ is a group algebra). This problem leads to
what we called the {\it Fourier analysis} of inductive families of
algebras (groups).

  We want to emphasize that the notion of $GZ_n$-subalgebra
of an algebra $A_n$ does depend on the structure of the inductive family
 $A_i$, $i=1,\dots, n$, and not only on the algebra $A_n$ itself; so if we choose
 another inductive family inside $A_n$, then $GZ_n$ also can change.
 The development of these ideas for the symmetric groups can be found
 in \cite{OV,VO}. The assumption on the simplicity of the spectrum is supposed
 to be satisfied in all further considerations.

The analysis of examples of Gelfand--Tsetlin algebras in the
case of groups, and especially of the $GZ_n$ subalgebras of
$C(S_N)$, allows us to formulate the following theorem.

\begin{theorem}
{Suppose that $G_1 \subset G_2\subset \dots \subset G_n $ is a finite
sequence of finite groups. Suppose that the restriction of irreducible
representations of $G_k$ to $G_{k-1}$, $k=1, \dots, n$, is
multiplicity-free
and there exists a virtual projection of $G_k$ to $G_{k-1}$, $k=1, \dots, n$.
Then
the family of sets $\{X_k =\ker P_k,\; k=1, \dots, n\}$
generates (as multiplicative generators) the subalgebra $GZ_n$; here $P_k$ is
the generalized expectation $C(G_k) \longrightarrow C(G_{k-1}$)
corresponding to the  virtual projection $p_k: G_k \to G_{k-1}$ (see
the previous section).}
\end{theorem}

\begin{proof}
{Using Remark 1 after Theorem 2, we can prove that the center of
$C(G_k)$ belongs to the algebra generated by $GZ_{k-1}$ and the
set $X_k$.}
\end{proof}

In the case of the symmetric group, the set $X_k$ is determined by
the YJM-elements.

\section{The inverse limit of an inductive family of algebras
and GZ-algebras,
and martingales}

Suppose now we have a countable sequence $A_n$, $n=0,1,
\dots$, $A_0=\bf C$, $A_n \subset A_{n+1}$, of  $C^*$-algebras that
form an inductive family of algebras and define the inductive
limit
  $$A_{\infty}=\limind A_i$$
with respect to the embedding of algebras.

  In the same spirit we can define the inductive limit of
  the Gelfand--Tsetlin algebras
 $$GZ_{\infty}=\limind GZ_n;$$
 under our assumptions, it
is again a maximal abelian
 subalgebra of $A_{\infty}$.

Using Theorem 4 from the previous section, we can define
multiplicative generators of $GZ_{\infty}=\mathop{\rm lim\,ind} GZ_n$ for
the case of group algebras. In particular, this gives a
description of a multiplicative basis for the $GZ$-algebra of the  infinite
symmetric group.

An inductive family $\{A_n\}$ of {\it finite-dimensional
algebras} determines a $\bf {Z}_+$-graded graph $Y$ (the Bratteli
diagram). The vertices of level $n\geq 0$ correspond to the
simple subalgebras of the algebra $A_n$ (at the zero level we have one
vertex $\bf 0$), and two adjacent levels $Y_n$ and $Y_{n-1}$
form precisely the bipartite graph that was mentioned in
Sec.~1. The set of all maximal paths  (finite if
the number of algebras is finite, or infinite) from the vertex $\bf
0$ to the end is called the set of {\it tableux} and denoted by $T(Y)$
(recall that a path is a sequence of edges, and in the
multiplicity-free case a path is also a sequence of vertices).
Now let us choose a sequence of generalized expectations
  of these algebras at each level:
  $$P_n:A_n \longrightarrow  A_{n-1},\quad n=1,2, \dots. $$

\begin{lemma}
{The restriction of the generalized expectation $P_n$ to the
Gelfand--Tsetlin algebra $GZ_n$ sends it to $GZ_{n-1}$; thus this
restriction is an expectation of $GZ_n$ to $GZ_{n-1}$.}
\end{lemma}

\begin{proof}
{Each expectation sends the center of the algebra onto the center of
the subalgebra: $P_n(\zeta(A_n)) = \zeta(A_{n-1})$. Indeed, let $z \in
\zeta(A_n)$ and $b \in A_{n-1}$; then
$P_n(zb)=P_n(z)b=P_n(bz)=bP_nz$. At the same time
$P_n(\zeta(A_{n-1}))=\zeta(A_{n-1})$. Consequently,
$P_n(GZ_n)=GZ_{n-1}$ by definition.}
\end{proof}

  Now let us define the projective limit
  $$A^{\infty}=\limproj\{A_n,P_n\}$$
 with respect to the sequence of generalized expectations.
 It is obvious from definition that the following
lemma holds.

\begin{lemma}
  $A^{\infty}$ is a left and right
  $A_{\infty}$-bimodule (but not an algebra in general).
\end{lemma}

 Indeed, all algebras $A_n$ act from the left and from the right on
 all $A_m$, $m>n$; thus these actions extend to the projective limit.
 This definition makes sense for a general inductive family
 with an arbitrary system of expectations.

By Lemma 3, we can also correctly define the inverse (projective) limit
of the algebras $\{GZ_n\}$:
   $$GZ^{\infty}=\limproj \{GZ_n,P_n\}.$$
This is not an algebra either,  but a module over $GZ_{\infty}$. The
interpretation of this limit will be given below.

 Suppose now that all algebras $A_n$, $n=1,2, \dots$,  are
finite-di\-men\-si\-on\-al semisimple algebras. Since (generalized)
expectations are determined by systems of (generalized) Markovian matrices,
the projective module is determined by the system of matrices
$\Lambda_n$, $n=1,2, \dots$, where $\Lambda_n$  determines the
expectation of $A_n$ to $A_{n-1}$.

Let us fix such a system of generalized (or true) Markovian
matrices $\Lambda_n$, $n=1,2, \dots $. The size of $\Lambda_n$ is
$m_n \times m_{n-1}$, where $m_k$ is the number of simple
subalgebras in the algebra $A_k$. We denote this system of matrices
by ${\bf L}=\{\Lambda_n,\; n=1, 2,\dots\}$, and in order to emphasize the
dependence of the projective limit on the expectations, we will
sometimes write
$$A^{\infty}_{\bf L}=\limproj \{A_n,P_n\}$$
and
$$GZ^{\infty}_{\bf L}=\limproj \{A_n,P_n\}.$$

In the case of abelian algebras, as well as in the case of
GZ-algebras, such an inverse limit is well known by another name,
at least when all matrices $\Lambda_n$ are true Markovian
matrices. We will shortly explain this link. First of all, as
usual, the system of Markovian matrices $\bf L$ determines a
Markov measure $\mu_{\bf L}$ on the space of tableaux $T(Y)$ (see
above). Thus we have a measure space (ore precisely, a Lebesgue
space) $(T(Y), \textbf{A} \mu_{\bf L})$, where $\textbf{A}$ is the
sigma-field generated by elementary cylindric sets (elementary
cylindric set of order $n$ is the set of all paths with common
fragment of length $n$). Second, in \textbf{A} we have an
increasing sequence of finite sigma-subfields of cylindric sets of
order $n$. Following the general definition of {\it martingales},
we can now define the vector space ${\bf M}_{\bf L}$  of
martingales over this increasing sequence of sigma-subfields, each
of them being a sequence $\{f_n\}_n$ of measurable functions such
that $f_n$ is $\textbf{A}_n$-measurable and the expectation of
$f_n$ on the sigma-field $\textbf{A}_{n-1}$ is equal to $f_{n-1}$.

It is clear from definition that {\it this space of
martingales is exactly the inverse limit
$GZ^{\infty}_{\bf L}$} defined above.

This is the reason for calling the elements of the inverse limit
$A^{\infty}_{\bf L}$ of algebras noncommutative martingales. This
opens a wide range of generalizations of the martingale theory to
this noncommutative case. If we have a generalized expectation,
then we need to consider martingales with respect to non-positive
measures, which, as far as I know, were never considered.

In the group case there is a distinguished Markov measure ---
the so-called Plancherel measure on the space of tableaux $T(Y)$;
namely, if $G= \mathop{\rm lim\,ind} G_n$ is a locally finite group with simple
spectrum (like $S_\infty=\mathop{\rm lim\,ind} S_n$), then, using one of the
expectations defined in the previous section, we obtain the
Plancherel measure on $T(Y)$, which is the inverse limit of the Plancherel
measures on the spaces of finite tableaux. Martingales with
respect to the Plancherel
measure play an important role as a special kind of modules over the
group algebras of the group $G$.

Our last remark concerns the link with von Neumann's theory of
infinite tensor products: if our algebra $A_{\infty}$ is the {\it
infinite tensor product of algebras of matrices} (e.g., of order
two), the so-called Glimm algebras, then each incomplete tensor
product of Hilbert spaces in the sense of \cite{vN} is generated
by the inverse limit of algebras with respect to some sequence of
expectations. In this spirit, the scheme of this section allows us
to generalize von Neumann's theory to an arbitrary inductive limit
of finite-dimensional algebras instead of Glimm algebras.

\end{document}